\documentclass{article}

\title{Graph components and dynamics over finite fields}
\author{Ryan Flynn\and Derek Garton}

\usepackage{amsfonts}
\usepackage{amsmath}
\usepackage{amsthm}
\usepackage{amssymb}
\usepackage{url}
\usepackage{hyperref}

\DeclareMathOperator{\PGL}{PGL}

\theoremstyle{remark}
\newtheorem{obs}{Observation}[section]
\newtheorem{Rats}[obs]{Remark}

\theoremstyle{plain}
\newtheorem{Quad}[obs]{Fact}
\newtheorem{Polys}[obs]{Lemma}
\newtheorem{PCycleBounds}[obs]{Corollary}
\newtheorem{PLowerBound}[obs]{Corollary}
\newtheorem{PUpperBound}[obs]{Corollary}
\newtheorem{PPeriods}[obs]{Corollary}

\newtheorem{RatBounds}[obs]{Lemma}
\newtheorem{Prov}[obs]{Lemma}
\newtheorem{PreCycleBounds}[obs]{Lemma}
\newtheorem{CycleBounds}[obs]{Corollary}
\newtheorem{LowerBound}[obs]{Theorem}
\newtheorem{UpperBound}[obs]{Corollary}
\newtheorem{Periods}[obs]{Corollary}

\begin{document}

\maketitle

\begin{abstract}
For polynomials and rational maps of fixed degree over a finite field, we bound both the average number of connected components of their functional graphs as well as the average number of periodic points of their associated dynamical systems.
\end{abstract}

\section{Introduction} \label{Introduction}

A \emph{(discrete) dynamical system} is a pair $\left(S,f\right)$ consisting of a set $S$ and a map $f:S\to S$. The \emph{functional graph} of such a system, $\Gamma(S,f)$, is the directed graph whose set of vertices is $S$ and whose edges are given by the relation $\alpha\to\beta$ if and only if $f(\alpha)=\beta$. A \emph{connected component} of such a graph is a connected component of the underlying undirected graph. If $S$ is finite, say $\vert S\vert=n$, then the average number of components of random functional graph on $S$ is
\[
n^{-n}\sum_{f:S\to S}{\left\vert\left\{\textrm{components of }\Gamma(S,f)\right\}\right\vert}.
\]
Kruskal (see~\cite{MR0062973}) proved that as $n\to\infty$, this quantity is simply
\[
\frac{1}{2}\log{n}+\left(\frac{\log{2}+C}{2}\right)+o(1),
\]
where $C=.5772\ldots$ is Euler's constant.

However, far less is known about the analogous situation for polynomials and rational functions over finite fields. More precisely, let $q=p^n$ for $n\in\mathbb{Z}^{>0}$ and define $\Gamma(q,f)=\Gamma\left(\mathbb{F}_q,f\right)$ if $f\in\mathbb{F}_q[x]$ and $\Gamma(q,f)=\Gamma\left(\mathbb{P}^1(\mathbb{F}_q),f\right)$ if $f\in\mathbb{F}_q(x)$. (If there is no ambiguity, we will frequently write $\Gamma_f$ for $\Gamma(q,f)$.) Then we can ask the question: what is the average number of connected components of $\Gamma_f$, for $f$ ranging over all polynomials (or rational functions) over $\mathbb{F}_q$ of a fixed degree? In particular, do these functional graphs behave like random functional graphs? In this paper, we obtain lower bounds (Theorems~\ref{PLowerBound} and~\ref{LowerBound}) and upper bounds (Corollaries~\ref{PUpperBound} and~\ref{UpperBound}) for the average number of connected components of such functional graphs. As corollaries of the theorems, we also find lower bounds on the average number of periodic points of such functional graphs (Corollaries~\ref{PPeriods} and~\ref{Periods}).

Our lower bounds, Theorems~\ref{PLowerBound} and~\ref{LowerBound}, show that as long as $d$ is large enough relative to $q$, the average number of connected components of graphs of polynomials (or rational functions) of degree $d$ over $\mathbb{F}_q$ is at least as great as the average number of connected components of random functional graphs on $q$ vertices as $q\to\infty$. Specifically, when $d\geq\sqrt{q}$, the average number of connected components of functional graphs of polynomials (or rational functions) of degree $d$ over $\mathbb{F}_q$ is bounded below by
\[
\frac{1}{2}\log{q}-4.
\]

Similarly, Corollaries~\ref{PPeriods} and~\ref{Periods} imply that when $d\geq\sqrt{q}$, the average number of periodic points of a polynomial (or rational function) of degree $d$ over $\mathbb{F}_q$ is bounded below by $\frac{5}{6}\sqrt{q}-4$. Harris (see~\cite{MR0119227}) proved the analog of Kruskal's result: the average number of periodic elements of a random functional graph on $n$ vertices is about
\[
\sqrt{\frac{\pi}{2}}\sqrt{n}
\]
as $n\to\infty$.

Before we begin our analysis, we need one more piece of notation. If $\alpha\in\mathbb{F}_q$ is periodic with orbit $\alpha=\alpha_1\to\alpha_2\to\cdots\to\alpha_k\to\alpha_1$ (with $\alpha_i$'s distinct), then we refer to its orbit as a \emph{cycle} (cycles of length $k$ are called \emph{$k$-cycles}). (See~\cite{MR2033734} for more exposition and illustrations of the cycle structure of functional graphs.) The basis for our results is the following observation:

\begin{obs} \label{obs}
The connected components of $\Gamma_f$ are in one-to-one correspondence with the cycles of $f$.
\end{obs}

Our method is simple: count the number of polynomial and rational functions of a fixed degree that give an arbitrary cycle---then sum over possible cycles to obtain our results. More precisely, if $d\in\mathbb{Z}^{\geq0}$, we compute the following quantities:
\[
\mathcal{P}(q,d):=\frac{\displaystyle{\sum_{\substack{f\in\mathbb{F}_q[x]\\\deg{\left(f\right)}=d}}}{\left\vert\left\{\textrm{cycles in }\Gamma_f\right\}\right\vert}}{\displaystyle{\sum_{\substack{f\in\mathbb{F}_q[x]\\\deg{\left(f\right)}=d}}}{1}},
\]
and, for any $k\in\mathbb{Z}^{>0}$,
\[
\mathcal{P}(q,d,k):=\frac{\displaystyle{\sum_{\substack{f\in\mathbb{F}_q[x]\\\deg{\left(f\right)}=d}}}{\left\vert\left\{k\textrm{-cycles in }\Gamma_f\right\}\right\vert}}{\displaystyle{\sum_{\substack{f\in\mathbb{F}_q[x]\\\deg{\left(f\right)}=d}}}{1}}.
\]
We define analogous quantities $\mathcal{R}(q,d)$ and $\mathcal{R}(q,d,k)$ for rational functions. We address the polynomial case in \hyperref[Polynomials]{Section~\ref*{Polynomials}} and the rational function case in \hyperref[Rationals]{Section~\ref*{Rationals}}. We note that our methods allow us to count only the ``short" cycles, so for our lower bounds we simply ignore the long ones. As it turns out, the long cycles are scarce enough that despite their omission, our component count still looks ``random".

Vasiga and Shallit~\cite{MR2033734} have previously studied the cycle structure of $\Gamma_f$ for the cases $f=x^2$ and $f=x^2-2$, as has Rogers~\cite{MR1368298} for $f=x^2$. The squaring function is also defined over $\mathbb{Z}/m\mathbb{Z}$; Carlip and Mincheva~\cite{MR2402530} have addressed this situation for certain $m$. Similarly, Chou and Shparlinski~\cite{MR2072394} have studied the cycle structure of repeated exponentiation over $\mathbb{F}_p$.

The question of whether polynomials act like random functions arises in the context of Pollard's rho algorithm for factoring integers (see~\cite{MR0392798}). To factor an integer $n$, the algorithm requires a pseudorandom function modulo prime divisors of $n$. Specifically, for a divisor $p\mid n$, a point should traverse a cycle after $O\left(\sqrt{p}\right)$ iterations of the function (this ``collision" produces a factor of $n$ divisible by $p$). In practice, it turns out that quadratic polynomials work for this purpose! That is, polynomials seem to have this property as soon as they have degree $d=2$. It is clear that polynomials over $\mathbb{F}_q$ of degree $d\geq q$ act like random functions on $\mathbb{F}_q$ for any reasonable definition of ``act like random functions", since such polynomials give every possible map on $\mathbb{F}_q$ the same number of times. In this paper, we show that polynomials (and rational functions) have as many cycles and periodic points as random functions when $d\geq\sqrt{q}$. We remark that the average time to collision of random functional graphs on $n$ vertices is $O(\sqrt{n})$, which follows from the fact that there are very few ``long cycles" in such graphs on average. The techniques of this paper, however, are not suited to controlling ``long cycles" (as will become apparent), so they do not lead to a useful ``collision time" result.

There are extensive data and heuristic arguments supporting the claim that degree-2 polynomials produce as many collisions as random functions, but very little has been proven (see~\cite{MR0392798} and~\cite{MR1094034}). (However, see~\cite{MR1094034} for a proof that a the expected number of steps to produce a factor of $n$ divisible $p$ is $\Omega\left(p^{-1}\log^2{p}\right)$.) One heuristic that suggests that degree-2 polynomials (or rational functions) act like random polynomials is the following fact.

\begin{Quad} \label{Quad}
Define a ``quadratic graph" to be any functional graph in which every vertex has either zero or two preimages. Then the average number of periodic points of a random quadratic graph on $n$ vertices is the same as that of a random graph, as $n\to\infty$.
\begin{proof}
Let $\Gamma$ be a directed graph on $n$ vertices such that the out-degree of each vertex is 1 and the in-degree is either 0 or $m$. Letting $n=mt$, there are $t$ vertices with in-degree $m$ and $(m-1)t$ vertices with in-degree 0. Assuming the vertices are labeled, the number of such graphs is
\[
\binom{mt}{t}\frac{(mt)!}{(m!)^t}.
\]
For $1\leq k\leq t$, the number of such graphs with a fixed $k$-cycle is
\[
\binom{mt-k}{t-k}\frac{(mt-k)!}{((m-1)!)^k(m!)^{t-k}},
\]
so the number of $k$-periodic points amongst all such graphs is
\[
k\cdot\binom{mt-k}{t-k}\frac{(mt-k)!}{((m-1)!)^k(m!)^{t-k}}\cdot\frac{(mt)!}{k(mt-k)!}=\binom{mt-k}{t-k}\frac{(mt)!}{((m-1)!)^k(m!)^{t-k}}.
\]
Summing over $1\leq k\leq t$, the average we want is
\[
\binom{mt}{t}^{-1}\sum_{k=1}^{t}{m^k\binom{mt-k}{t-k}}=-1+\binom{mt}{t}^{-1}\sum_{k=0}^{t}{m^k\binom{mt-k}{t-k}}.
\]
When $m=1$, the average is $t$ (since all points would be periodic). When $m=2$, as in the question, the sum is 
\[
-1+\frac{4^t}{\binom{mt}{t}}\sim-1+\sqrt{\frac{\pi}{2}}\sqrt{n},
\]
as desired.
\end{proof}
\end{Quad}

\section{Polynomials} \label{Polynomials}

\begin{Polys} \label{Polys}
If $1\leq k\leq d+1$, then
\[
\sum_{\substack{f\in\mathbb{F}_q[x]\\\deg{\left(f\right)}\leq d}}{\left\vert\left\{k\textrm{-cycles in }\Gamma_f\right\}\right\vert}
=\frac{q\left(q-1\right)\cdots\left(q-(k-1)\right)}{k}\cdot q^{d-k+1}.
\]
\begin{proof}
If $k>q$, then the statement is trivially true. Thus, suppose that $k\leq q$. Now, fix any $k$-cycle $C:\alpha_1\to\cdots\to\alpha_k\to\alpha_1$. Note that $f\in\mathbb{F}_q[x]$ gives $C$ if and only if $f(\alpha_i)=\alpha_{i+1}$ for all $1\leq i\leq k-1$ and $f(\alpha_k)=\alpha_1$. Lagrange interpolation produces such a polynomial of degree $k-1$ or less; call it $f_C$. Moreover, if $g\in\mathbb{F}_q[x]$ vanishes at every $\alpha_i$, then $f_C+g$ gives $C$ as well. Additionally, every polynomial of degree $d$ or less giving $C$ is of this form for some $g\in\mathbb{F}_q[x]$ of degree $d$ or less. By linear algebra, the number of polynomials of degree $d$ or less in $\mathbb{F}_q[x]$ that vanish at every $\alpha_i$ is $q^{d+1-k}$; since there are $k^{-1}q\left(q-1\right)\cdots\left(q-k+1\right)$ $k$-cycles, the lemma is proved.
\end{proof}
\end{Polys}

\begin{PCycleBounds} \label{PCycleBounds}
If $1\leq k\leq d$ or if $0=d=k-1$, then
\[
\mathcal{P}(q,d,k)=\frac{q(q-1)\cdots(q-(k-1))}{kq^k}.
\]

\begin{proof}
Note that
\[
\mathcal{P}(q,d,k)=
\frac{\displaystyle{\sum_{\substack{f\in\mathbb{F}_q[x]\\\deg{\left(f\right)}\leq d}}{\left\vert\left\{k\textrm{-cycles in }\Gamma_f\right\}\right\vert}
-\sum_{\substack{f\in\mathbb{F}_q[x]\\\deg{\left(f\right)}\leq d-1}}{\left\vert\left\{k\textrm{-cycles in }\Gamma_f\right\}\right\vert}}}
{q^{d+1}-q^d},
\]
then apply \hyperref[Polys]{Lemma~\ref*{Polys}}.
\end{proof}
\end{PCycleBounds}

\begin{PLowerBound} \label{PLowerBound}
\[
\mathcal{P}(q,d)\geq\sum_{k=1}^{\min{\{d,q\}}}{\frac{q(q-1)\cdots(q-(k-1))}{kq^k}},
\]
with equality if and only if $d\geq q$. In particular,
\begin{align*}
\mathcal{P}(q,d)&>\sum_{k=1}^{\min{\{d,\lfloor\sqrt{q}\rfloor\}}}{\left(\frac{1}{k}\right)}-\frac{1}{4}\\
&>\log{\left(\min{\{d,\lfloor\sqrt{q}\rfloor\}}+1\right)}-\frac{1}{4}.
\end{align*}
\begin{proof}
Immediate by \hyperref[PCycleBounds]{Corollary~\ref*{PCycleBounds}}.
\end{proof}
\end{PLowerBound}

\begin{PUpperBound} \label{PUpperBound}
If $d<q-1$, then
\[
\mathcal{P}(q,d)\leq\sum_{k=1}^{d+1}{\left(\frac{1}{k}\right)}+\frac{q}{d+2}.
\]
If $d\geq q-1$, then
\[
\mathcal{P}(q,d)\leq\sum_{k=1}^{q}{\left(\frac{1}{k}\right)}.
\]

\begin{proof}
Using \hyperref[PCycleBounds]{Corollary~\ref*{PCycleBounds}} and \hyperref[Polys]{Lemma~\ref*{Polys}}, we know that
\begin{align*}
\mathcal{P}(q,d)&=\sum_{k=1}^{q}{\mathcal{P}(q,d,k)}\\
&\leq\sum_{k=1}^{d}{\left(\frac{1}{k}\right)}+\mathcal{P}(q,d,d+1)+\sum_{k=d+2}^{q}{\mathcal{P}(q,d,k)}\\
&<\sum_{k=1}^{d}{\left(\frac{1}{k}\right)}+\frac{1}{d+1}+\sum_{k=d+2}^{q}{\mathcal{P}(q,d,k)}
\end{align*}
so it remains only to bound the frequency that ``long'' cycles appear (if $q\leq d+1$, there are no ``long'' cycles, so the second statement of the proposition is proved).

But a polynomial can have most $q(d+2)^{-1}$ cycles of length at least $(d+2)$, so the the proposition is proved.
\end{proof}
\end{PUpperBound}

\begin{PPeriods} \label{PPeriods}
The average number of periodic points of degree $d$ polynomials in $\mathbb{F}_q[x]$ is at least
\[
\sum_{k=1}^{\min{\{d,q\}}}{\frac{q(q-1)\cdots(q-(k-1))}{q^k}},
\]
with equality if and only if $d\geq q$. In particular, the average number of periodic points of degree $d$ polynomials in $\mathbb{F}_q[x]$ is at least $\frac{5}{6}\min{\{d,\lfloor\sqrt{q}\rfloor\}}$.

\begin{proof}
Note that a $k$-cycle contains $k$ periodic points and that
\[
\sum_{k=1}^{\lfloor\sqrt{q}\rfloor}{\left(1-\frac{k^2-k}{2q}\right)}=\lfloor\sqrt{q}\rfloor-\frac{\lfloor\sqrt{q}\rfloor^3-\lfloor\sqrt{q}\rfloor}{6q}>\frac{5}{6}\lfloor\sqrt{q}\rfloor,
\]
then use \hyperref[PCycleBounds]{Corollary~\ref*{PCycleBounds}}.
\end{proof}
\end{PPeriods}

\section{Rational Functions} \label{Rationals}

\begin{Rats} \label{Rats}
We first make a distinction between rational functions on $\mathbb{P}^1(\mathbb{F}_q)$ and elements  of $\mathbb{F}_q(x)$: the first set contains the constant function sending $\mathbb{P}^1(\mathbb{F}_q)$ to $\infty$ whereas the second set does not. From the dynamical perspective, the constant infinity function is a perfectly valid rational map, so we will include it in our counts. Moreover, we say that it, like all other constant functions, has degree 0.
\end{Rats}

\begin{RatBounds} \label{RatBounds}
The number of rational functions on $\mathbb{P}^1(\mathbb{F}_q)$ of degree at most $d$ is $q^{2d+1}+1$. In particular, the number of rational functions on $\mathbb{P}^1(\mathbb{F}_q)$ of degree exactly $d$ is either $q+1$ (when $d=0$) or $q^{2d-1}\left(q^2-1\right)$ (when $d>0$).
\begin{proof}
There are $q^{2d+2}$ possible pairs of polynomials of degree at most $d$. The probability that these polynomials are relatively prime is $1-\frac{1}{q}+\frac{q-1}{q^{2d+2}}$ (see Corollary~4 of~\cite{MR2322084}, for example). To complete the count, we must divide out by scalars. Thus, the number of rational functions of degree at most $d$ is
\[
\frac{q^{2d+2}\left(1-\frac{1}{q}+\frac{q-1}{q^{2d+2}}\right)}{q-1}=q^{2d+1}+1,
\]
as desired.
\end{proof}
\end{RatBounds}

We will need a tool for counting rational functions taking specified values on points of $\mathbb{P}^1(\mathbb{F}_q)$; the main technical tool to this end is the following lemma. We remark that the ``long cycle" phenomenon impedes our counting: that is the reason we obtain an inequality when $j<m$ in the following lemma, and that is the reason we restrict degrees in the definitions of $S(d,C)^{\textrm{red}}$ and $S(d,C,g_0)$ in \hyperref[PreCycleBounds]{Lemma~\ref*{PreCycleBounds}}. In these instances, the ``long cycle" phenomenon manifests itself when counting rational functions of low degree: when the cycle $C$ is too long compared to the degree of a a denominator, it is difficult to say whether there exists a numerator of the same degree such that their quotient gives $C$.

\begin{Prov} \label{Prov}
Choose $\beta_1,\ldots,\beta_m\in\mathbb{F}_q$ distinct and choose any $\gamma_1,\ldots,\gamma_m\in\mathbb{F}_q$. Choose a monic $G_0\in\mathbb{F}_q[x]$ of degree $j_0$ that is either constant or irreducible, and choose a monic $G_1\in\mathbb{F}_q[x]$ of degree $j_1$. Let $G=G_0G_1$ and $j=\deg{(G)}=j_0+j_1$. If
\[
S\left(G_0,G_1\right)=\left\{f\in\mathbb{F}_q[x]\,\Big\vert\,
\substack{f\textrm{ monic},\,\deg{(f)}=j,\,G_0\mid f\\
f(\beta_i)=\gamma_iG(\beta_i)\textrm{ for all }i\in\{1,\ldots,m\}}\right\},
\]
then
\[
\left\vert S\left(G_0,G_1\right)\right\vert\leq
\begin{cases}
1&\textrm{if }j<m\\
1&\textrm{if }j_1<m\leq j\textrm{ and }G_0\neq\left(x-\beta_i\right)\textrm{ for any }i\in\{1,\ldots,m\}\\
q^{j-m}&\textrm{if }j\geq m\textrm{ and either }G_0=1\textrm{ or }G_0=\left(x-\beta_i\right)\textrm{ for some }i\in\{1,\ldots,m\}\\
q^{j_1-m}&\textrm{if }j_1\geq m\textrm{ and }G_0\neq\left(x-\beta_i\right)\textrm{ for any }i\in\{1,\ldots,m\},
\end{cases}
\]
with equality in the latter two cases.
\begin{proof}
As in \hyperref[Polynomials]{Section~\ref*{Polynomials}}, there exists a unique polynomial $f_C$ of degree less than $m$ such that $f_C(\beta_i)=\gamma_iG(\beta_i)$ for all $i\in\{1,\ldots,m\}$. If $j<m$, $f_C$ may or may not be monic, $f_C$ may or may not have degree $j$, and $G_0$ may or may not divide $f_C$; this gives the first case. However, if $j\geq m$, then
\[
\left\vert\left\{f\in\mathbb{F}_q[x]\,\Big\vert\,
\substack{f\textrm{ monic},\,\deg{(f)}=j,\\
f(\beta_i)=\gamma_iG(\beta_i)\textrm{ for all }i\in\{1,\ldots,m\}}\right\}\right\vert=
\left\vert\left\{h\in\mathbb{F}_q[x]\,\Big\vert\,
\substack{h\textrm{ monic},\\
\deg{(h)}=j-m}\right\}\right\vert=
q^{j-m},
\]
since those $f\in\mathbb{F}_q[x]$ of degree $j$ such that $f_C(\beta_i)=\gamma_iG(\beta_i)$ for all $i\in\{1,\ldots,m\}$ are exactly those of the form $f=f_h=f_C+h\prod_{i=1}^m{\left(x-\beta_i\right)}$ for some monic $h$ of degree $j-m$. We wil count those $f_h$ that are divisible by $G_0$.

First, note that if either $G_0=1$ or $G_0=\left(x-\beta_i\right)$ for any $i\in\{1,\ldots,m\}$, then $G_0$ divides every $f_h$, giving the third case; thus, suppose that $G_0$ is not of this form. Next, let $f_0\in\mathbb{F}_q[x]$ be the unique polynomial such that $f_0$ reduces modulo $G_0$ to $-f_C\prod_{i=1}^m{\left(x-\beta_i\right)^{-1}}$ and such that $\deg{\left(f_0\right)}\leq j_0-1$. Then $G_0\mid f_h$ if and only if $h$ is congruent to $f_0$ modulo $G_0$; i.e., if and only if $h=f_0+bG_0$ for some polynomial $b$. If $j_1<m$, then $\deg{\left(h\right)}=j-m<j_0$, so the only possibility for $b$ is $b=0$, giving the second case. If, on the other hand, $j_1\geq m$, then $\deg{\left(h\right)}\geq j_0$, so these $b$ are exactly the monic polynomials of degree $j-m-j_0=j_1-m$. Thus, there are $q^{j_1-m}$ of them, giving the final case.
\end{proof}
\end{Prov}

We can now turn to the question of bounds for the number of rational functions giving $k$-cycles for $1\leq k\leq d+1$; that is, the rational function analog of \hyperref[Polys]{Lemma~\ref*{Polys}}. The following lemma allows us to show (in \hyperref[CycleBounds]{Corollary~\ref*{CycleBounds}}) that $\mathcal{R}(q,d,k)$ is what we expect: about $k^{-1}$, at least for $k$ small compared to $q$ (just as in the polynomial case).

\begin{PreCycleBounds} \label{PreCycleBounds}
Let $K_{q,d,k}=\frac{(q+1)q\cdots(q+1-(k-1))}{k}q^{2d-k}$. If $1\leq k\leq d+1$, then
\[
(q-k-1)K_{q,d,k}<
\sum_{\substack{f\in\mathbb{F}_q(x)\\\deg{\left(f\right)}\leq d}}{\left\vert\left\{k\textrm{-cycles in }\Gamma_f\right\}\right\vert}
<qK_{q,d,k}.
\]
\begin{proof}

We will count how many rational functions give a fixed $k$-cycle $C:\alpha_1\to\cdots\to\alpha_k\to\alpha_1$. The result will then follow from the fact that there are $k^{-1}(q+1)\cdot q\cdots(q+1-(k-1))$ possible $k$-cycles.
Note that the group $\PGL_2{\left(\mathbb{F}_q\right)}$ acts via postcomposition on the set of rational functions on $\mathbb{P}^1\left(\mathbb{F}_q\right)$. For $k=1,2,3$ it is clear that exactly $\frac{1}{q+1},\frac{1}{q(q+1)},\frac{1}{q(q^2-1)}$ of the functions in every orbit give $C$ (respectively), so the result follows. Thus, we assume that $k\geq4$.
Moreover, we assume that $\alpha_k=\infty$ and $\alpha_1=1$; this restriction is harmless, because $f$ gives $\phi(C)$ if and only if $f^{\phi}$ gives $C$ (where $f^{\phi}:=\phi\circ f\circ\phi^{-1}$ is the conjugate of $f$ by $\phi$).

We will estimate the number of rational functions that give $C$ by summing over possible monic denominators. Since $\alpha_{k}=\infty$, such denominators must be multiples of $\left(x-\alpha_{k-1}\right)$. Thus, the set whose size we must estimate is:
\[
S(d,C):=\left\{(f,g)\in\left(\mathbb{F}_q[x]\right)^2\,\Big\vert\,
\substack{f\textrm{ and }g\textrm{ monic},\,\deg{(f)}=\deg{(g)}\leq d,\,\gcd{(f,g)}=1,\\
\,f(\alpha_i)=\alpha_{i+1}g(\alpha_i)\textrm{ for all }i\in\{1,\ldots,k-2\}}\right\}.
\]
Recall that our counting lemma, \hyperref[Prov]{Lemma~\ref*{Prov}}, does not provide an exact answer when the degree of the denominator is small relative to the length of the cycle. We can avoid this problem when bounding $\left\vert S(d,C)\right\vert$ from below, however, by simply omitting such low-degree denominators from our count. To this end, let
\[
S(d,C)^{\textrm{red}}=\left\{(f,g_1)\in\left(\mathbb{F}_q[x]\right)^2\,\Big\vert\,
\substack{f\textrm{ and }g_1\textrm{ monic},\,k-2\leq\deg{(f)}=1+\deg{\left(g_1\right)}\leq d,\,(x-\alpha_{k-1})\nmid f\\
f(\alpha_i)=\alpha_{i+1}(\alpha_i-\alpha_{k-1})g_1(\alpha_i)\textrm{ for all }i\in\{1,\ldots,k-2\}}\right\},
\]
which contains pairs of the form $(f,(x-\alpha_{k-1})g_1)$, whose components may or may not have a factor in common and such that the degrees of $f$ and $(x-\alpha_{k-1})g_1$ are at least $k-2$. To address the possibility of common factors, let
\[
S(d,C,g_0)=\left\{(f,g_1)\in\left(\mathbb{F}_q[x]\right)^2\,\Big\vert\,
\substack{f\textrm{ and }g_1\textrm{ monic},\,k-2\leq\deg{(f)}=j_0+\left(1+\deg{(g_1)}\right)\leq d,g_0\mid f\\
f(\alpha_i)=\alpha_{i+1}g_0\left(\alpha_i\right)(\alpha_i-\alpha_{k-1})g_1(\alpha_i)\textrm{ for all }i\in\{1,\ldots,k-2\}}\right\},
\]
where $g_0$ is monic and irreducible with $\deg(g_0)=j_0>1$; this set contains pairs of the form $(g_0f, g_0(x-\alpha_{k-1})g_1)$, that is, the $(f,g)\in S(d,C)^{\textrm{red}}$ such that $\gcd{(f,g)}$ is a multiple of $g_0$. We can estimate the size of both the above sets using \hyperref[Prov]{Lemma~\ref*{Prov}}; using that lemma's notation, we let $S_C(G_0,G_1)$ be associated to the data $m=k-2$, $\beta_i=\alpha_i$, and $\gamma_i=\alpha_{i+1}$ for $i\in\{1,\ldots,m\}$.

To bound $\left\vert S(d,C)\right\vert$ from below, we will estimate
\[
\left\vert S(d,C)^{\textrm{red}}\right\vert-\sum_{g_0\textrm{ irreducible}}{\left\vert S(d,C,g_0)\right\vert}.
\]
To estimate $\left\vert S(d,C)^{\textrm{red}}\right\vert$, note that for any monic $g_1$ with $m\leq\deg{\left(g_1\right)}=j_1\leq d-1$ we know by \hyperref[Prov]{Lemma~\ref*{Prov}} that
\[
\left\vert S_C\left(1,\left(x-\alpha_{k-1}\right)g_1\right)\right\vert
-\left\vert S_C\left(\left(x-\alpha_{k-1}\right),g_1\right)\right\vert
>q^{j_1+1-m}-q^{j_1-m},
\]
so summing over such $g_1$s gives
\begin{align*}
\left\vert S(d,C)^{\textrm{red}}\right\vert&>\sum_{j_1=m}^{d-1}{q^{j_1}\left(q^{j_1+1-m}-q^{j_1-m}\right)}\\
&=q^m\left(q-1\right)\left(\frac{q^{2(d-m)}-1}{q^2-1}\right).
\end{align*}
We must spit the estimation of $\left\vert S(d,C,g_0)\right\vert$ into two cases: whether $g_0$ is of the form $g_0=\left(x-\alpha_i\right)$ for some $i\in\{1,\ldots,m\}$ or not. If $g_0$ is of this form, then \hyperref[Prov]{Lemma~\ref*{Prov}} implies that if $g_1$ is monic with $m-2\leq\deg{\left(g_1\right)}=j_1\leq d-2$, then
\[
\left\vert S_C\left(g_0,\left(x-\alpha_{k-1}\right)g_1\right)\right\vert=q^{j_1+2-m},
\]
so
\[
\sum_{i=1}^{m}{\left\vert S_C\left(\left(x-\alpha_i\right),\left(x-\alpha_{k-1}\right)g_1\right)\right\vert}
=m\cdot\sum_{j_1=m-2}^{d-2}{q^{j_1}q^{j_1+2-m}}
=mq^{m-2}\left(\frac{q^{2(d-m+1)}-1}{q^2-1}\right)
\]
If, on the other hand, $g_0$ is not of this form and $\deg{\left(g_0\right)}=j_0$, then \hyperref[Prov]{Lemma~\ref*{Prov}} implies that if $g_1$ is monic with $m-1\leq\deg{\left(g_1\right)}=j_1\leq d-1-j_0$, then
\[
\left\vert S_C\left(g_0,\left(x-\alpha_{k-1}\right)g_1\right)\right\vert=q^{j_1+1-m}.
\]
If $j_1<m-1$, then
\[
\left\vert S_C\left(g_0,\left(x-\alpha_{k-1}\right)g_1\right)\right\vert\leq1,
\]
so
\begin{align*}
&\sum_{\substack{g_0\textrm{ irreducible}\\g_0\neq\left(x-\alpha_i\right)\forall i\in\left\{1,\ldots,m\right\}}}{\left\vert S(d,C,g_0)\right\vert}\\
&\qquad<\sum_{j_0=1}^{d-m}{q^{j_0}\left(\sum_{j_1=0}^{m-2}{q^{j_1}}+\sum_{j_1=m-1}^{d-1-j_0}{q^{j_1}q^{j_1+1-m}}\right)}+\sum_{j_0=d-m+1}^{d-1}{q^{j_0}\sum_{j_1=0}^{d-1-j_0}{q^{j_1}}}\\
&\qquad<q\cdot\frac{q^{d-m}-1}{q-1}\cdot\frac{q^{m-1}-1}{q-1}+\frac{q^{m+1}}{q-1}\cdot\frac{q^{2(d-m)}-1}{q^2-1}+(m-1)q^{d-m+1}\cdot\frac{q^{m-1}-1}{q-1}\\
&\qquad<\frac{q^{d+1}}{(q-1)(q^2-1)}\cdot\left(q^{d-m}+m+1\right).
\end{align*}
Thus, we see that
\begin{align*}
&\sum_{g_0\textrm{ irreducible}}{\left\vert S(d,C,g_0)\right\vert}\\
&\qquad<mq^{m-2}\left(\frac{q^{2(d-m+1)}-1}{q^2-1}\right)+\frac{q^{d+1}}{(q-1)(q^2-1)}\cdot\left(q^{d-m}+m+1\right)\\
&\qquad<(m+1)\cdot\frac{q^{2d-m+1}}{(q-1)\left(q^2-1\right)}.
\end{align*}

Finally, we obtain the lower bound of
\begin{align*}
\left\vert S(d,C)\right\vert&>q^m\left(q-1\right)\left(\frac{q^{2(d-m)}-1}{q^2-1}\right)
-(m+1)\cdot\frac{q^{2d-m+1}}{(q-1)\left(q^2-1\right)}\\
&>\frac{q^{2d-m+1}}{(q-1)\left(q^2-1\right)}\left(q-(m+3)\right)\\
&>q^{2d-k}\left(q-k-1\right).
\end{align*}

To bound $\left\vert S(d,C)\right\vert$ from above, we must address the low-degree rational functions we omitted for our lower bound (see definition of $S(d,C)^{\textrm{red}}$). By \hyperref[Prov]{Lemma~\ref*{Prov}} we know that if $g_1$ is monic with $0\leq\deg{\left(g_1\right)}=m-1$, then $\left\vert S_C\left(1,\left(x-\alpha_{k-1}\right)g_1\right)\right\vert\leq1.$ For higher degree $g_1$, we sum as in the lower bound, without worrying about relative primality of the numerator and denominators. We obtain an upper bound of
\begin{align*}
&\sum_{j_1=0}^{m-1}{q^{j_1}}+\sum_{j_1=m}^{d-1}{q^{j_1}\left(q^{j_1+1-m}-q^{j_1+m}\right)}\\
&\qquad=\frac{q^m-1}{q-1}+q^m\left(q-1\right)\left(\frac{q^{2(d-m)}-1}{q^2-1}\right)\\
&\qquad<q^{2d-k+1},
\end{align*}
as desired.
\end{proof}
\end{PreCycleBounds}

Recall that for any $k\in\mathbb{Z}^{>0}$,
\[
\mathcal{R}(q,d,k):=\frac{\displaystyle{\sum_{\substack{f\in\mathbb{F}_q(x)\\\deg{\left(f\right)}=d}}}{\left\vert\left\{k\textrm{-cycles in }\Gamma_f\right\}\right\vert}}{\displaystyle{\sum_{\substack{f\in\mathbb{F}_q(x)\\\deg{\left(f\right)}=d}}}{1}}.
\]

\begin{CycleBounds} \label{CycleBounds}
If $1\leq k\leq d$, then
\[
\mathcal{R}(q,d,k)>
\frac{(q+1)q\cdots(q+1-(k-1))}{kq^k}\cdot\left(1-\frac{2k+2}{q}\right)
\]
and
\[
\mathcal{R}(q,d,k)
<\frac{(q+1)q\cdots(q+1-(k-1))}{kq^k}\cdot\left(1+\frac{2}{q^2}\right).
\]
In fact, the upper bound holds for $k=d+1$.

\begin{proof}

For the lower bound, we use Lemmas~\ref{RatBounds} and~\ref{PreCycleBounds} to see that
\begin{align*}
\mathcal{R}(q,d,k)&=\frac{\displaystyle{\sum_{\substack{f\in\mathbb{F}_q(x)\\\deg{\left(f\right)}\leq d}}}{\left\vert\left\{k\textrm{-cycles in }\Gamma_f\right\}\right\vert-\sum_{\substack{f\in\mathbb{F}_q(x)\\\deg{\left(f\right)}\leq d-1}}}{\left\vert\left\{k\textrm{-cycles in }\Gamma_f\right\}\right\vert}}{q^{2d-1}\left(q^2-1\right)}\\
&>\frac{\prod_{j=0}^{k-1}{(q+1-j)}}{kq^k}
\cdot\frac{q^2-(k+1)q-1}{q^2-1}\\
&=\frac{\prod_{j=0}^{k-1}{(q+1-j)}}{kq^k}\cdot\left(1-(k+1)\frac{q}{q^2-1}\right),
\end{align*}
as desired.

For the upper bound, use Lemmas~\ref{RatBounds} and~\ref{PreCycleBounds} again to see that
\begin{align*}
\mathcal{R}(q,d,k)
&<\frac{\prod_{j=0}^{k-1}{(q+1-j)}}{kq^k}
\cdot\frac{q^2}{q^2-1}\\
&=\frac{\prod_{j=0}^{k-1}{(q+1-j)}}{kq^k}\cdot\left(1+\frac{1}{q^2-1}\right),
\end{align*}
which completes the proof.
\end{proof}
\end{CycleBounds}

\begin{LowerBound} \label{LowerBound}
\begin{align*}
\mathcal{R}(q,d)&>\sum_{k=1}^{\min{\{d,\lfloor\sqrt{q}\rfloor\}}}{\left(\frac{1}{k}\right)}-4\\
&>\log{\left(\min{\{d,\lfloor\sqrt{q}\rfloor\}}+1\right)}-4.
\end{align*}

\begin{proof}
Using \hyperref[CycleBounds]{Corollary~\ref*{CycleBounds}}, it follows that
\begin{align*}
\mathcal{R}(q,d)
&\geq\sum_{k=1}^{d}{\mathcal{R}(q,d,k)}\\
&>\sum_{k=1}^{d}{\frac{q(q-1)\cdots(q-(k-1))}{kq^k}\left(1-\frac{2k+2}{q}\right)}\\
&\geq\sum_{k=1}^{d}{\left(\frac{1}{k}-\frac{k-1}{2q}\right)\left(1-\frac{2k+2}{q}\right)}\\
&>\sum_{k=1}^{d}{\left(\frac{1}{k}-\frac{k-1}{2q}-\frac{2k+2}{kq}\right)}\\
&=\sum_{k=1}^{d}{\left(\frac{1}{k}-\frac{1}{2q}\left(k-1+\frac{4k+4}{k}\right)\right)}\\
&\geq\sum_{k=1}^{d}{\left(\frac{1}{k}\right)}-\frac{1}{2q}\sum_{k=1}^{d}{\left(k+7\right)}\\
&=\sum_{k=1}^{d}{\left(\frac{1}{k}\right)}-\frac{d^2+15d}{4q},
\end{align*}
giving the result.
\end{proof}
\end{LowerBound}

\begin{UpperBound} \label{UpperBound}
If $d<q-1$, then
\[
\mathcal{R}(q,d)\leq\frac{q+1}{d+2}+\left(1+\frac{2}{q^2}\right)\sum_{k=1}^{d+1}{\left(\frac{1}{k}\right)}.
\]
If $d\geq q-1$, then
\[
\mathcal{R}(q,d)\leq\left(1+\frac{2}{q^2}\right)\sum_{k=1}^{q+1}{\left(\frac{1}{k}\right)}.
\]

\begin{proof}
Using \hyperref[CycleBounds]{Corollary~\ref*{CycleBounds}} and \hyperref[PreCycleBounds]{Lemma~\ref*{PreCycleBounds}}, we see that
\begin{align*}
\mathcal{R}(q,d)&=\sum_{k=1}^{d+1}{\mathcal{R}(q,d,k)}
+\sum_{k=d+2}^{q+1}{\mathcal{R}(q,d,k)}\\
&<\left(1+\frac{2}{q^2}\right)\sum_{k=1}^{d+1}{\left(\frac{1}{k}\right)}
+\sum_{k=d+2}^{q}{\mathcal{P}(q,d,k)},
\end{align*}
so it remains only to bound the frequency with which ``long'' cycles appear (here, long cycles are those of length at least $d+2$).

As in \hyperref[PUpperBound]{Corollary~\ref*{PUpperBound}}, a rational function on $\mathbb{P}^1(\mathbb{F}_q)$ can have most $(q+1)(d+2)^{-1}$ cycles of length at least $(d+2)$, so the the proposition is proved.
\end{proof}
\end{UpperBound}

\begin{Periods} \label{Periods}
The average number of periodic points of degree $d$ rational functions on $\mathbb{P}^1\left(\mathbb{F}_q\right)$ is at least
\[
\sum_{k=1}^{\min{\{d,q\}}}{\left(\frac{(q+1)q\cdots(q-(k-1))}{q^k}\right)\left(1-\frac{k+4}{q}\right)},
\]
In particular, the average number of periodic points of degree $d$ rational functions on $\mathbb{P}^1\left(\mathbb{F}_q\right)$ is at least $\frac{5}{6}\min{\{d,\lfloor\sqrt{q}\rfloor\}}-3$.
\begin{proof}
Note that a $k$-cycle contains $k$ periodic points and that
\begin{align*}
\sum_{k=1}^{\lfloor\sqrt{q}\rfloor}{\left(1-\frac{k^2-k}{2q}\right)\left(1-\frac{k+2}{q}\right)}
&>\sum_{k=1}^{\lfloor\sqrt{q}\rfloor}{\left(1-\frac{k^2-k}{2q}-\frac{k+2}{q}\right)}\\
&=\lfloor\sqrt{q}\rfloor-\frac{2\lfloor\sqrt{q}\rfloor^3+6\lfloor\sqrt{q}\rfloor^2+28\lfloor\sqrt{q}\rfloor}{12q},
\end{align*}
then use \hyperref[CycleBounds]{Corollary~\ref*{CycleBounds}}.
\end{proof}
\end{Periods}

\section{Acknowledgements} \label{Acknowledgements}

We would like to thank the organizers of the 2010 Arizona Winter School on Number Theory and Dynamics (\url{http://math.arizona.edu/~swc/aws/10/}) for an educational week of arithmetic dynamics. Without the generous funding from the \href{http://swc.math.arizona.edu/}{AWS} and the help from the mathematicians we met there, we would not have written this paper. In particular, we would like to thank Professors \href{http://www.math.brown.edu/~jhs/}{Joseph Silverman}, \href{http://mathcs.holycross.edu/~rjones/}{Rafe Jones}, and \href{http://www.math.hawaii.edu/~mmanes/}{Michelle Manes} for suggesting this project and providing helpful advice while working on it.

The second author is also grateful for partial support of his work provided by the \href{http://www.nsf.gov/awardsearch/showAward.do?AwardNumber=0838210&version=noscript}{NSF-RTG} grant ``Algebraic Geometry and Number Theory at the University of Wisconsin".

\bibliography{graphcomponentsbib}
\bibliographystyle{amsplain}

\end{document}